\documentclass[11pt]{article}
\usepackage{bbm}
\usepackage{mathrsfs}
\usepackage{amscd}
\usepackage{amsmath,amsfonts,amssymb,amscd}
\usepackage{indentfirst,graphics,epsfig,psfrag}
\input{epsf}
\usepackage{ifpdf}
\usepackage{enumerate}
\usepackage{appendix}
\usepackage{enumerate}
\usepackage{amsmath}
\usepackage{lineno}
\usepackage{color}
\makeatletter \@addtoreset{figure}{section} \makeatother
\makeatletter
\long\def\@makecaption#1#2{%
   \vskip 10\p@
   \setbox\@tempboxa\hbox{{#1}\ \ #2}%
   \ifdim \wd\@tempboxa >\hsize

       {#1}\ \ #2\par
   \else
       \hbox to\hsize{\hfil\box\@tempboxa\hfil}%
   \fi}
\makeatother

\newtheorem{thm}{Theorem}
\newtheorem{cor}{Corollary}
\newtheorem{lem}{Lemma}

\newtheorem{pro}{Proposition}

\newcommand{\qed}{{\hfill\rule{3pt}{7pt}}}

\def\qed{\hfill \rule{4pt}{7pt}}

\begin{document}
\title{\textbf{On the equitable vertex arboricity of complete tripartite graphs}
\footnote{Research supported by the National Science Foundation of
China (Nos. 61164005 and 11161037), the National Basic Research
Program of China (No. 2010CB334708) and the Program for Changjiang
Scholars and Innovative Research Team in Universities (No. IRT1068),
the Research Fund for the Chunhui Program of Ministry of Education
of China (No. Z2014022) and the Nature Science Foundation from
Qinghai Province (Nos. 2012-Z-943, 2014-ZJ-721 and 2014-ZJ-907).}}
\author{
\small Zhiwei Guo$^{1}$ Haixing Zhao$^2$ Yaping Mao$^1$\\[0.3cm]
\small $^1$Department of Mathematics, Qinghai Normal\\
\small University, Xining, Qinghai 810008, China\\
\small E-mails: guozhiweic@yahoo.com; maoyaping@ymail.com\\[0.3cm]
\small $^2$School of Computer, Qinghai Normal\\
\small University, Xining, Qinghai 810008, China\\
\small E-mails: h.x.zhao@163.com \\[0.3cm]}
\date{}
\maketitle
\begin{abstract}
The equitable coloring problem, introduced by Meyer in 1973,
has received considerable attention and research. Recently, Wu et
al. introduced the concept of equitable $(t, k)$-tree-coloring,
which can be viewed as a generalization of proper equitable
$t$-coloring. The strong equitable
vertex $k$-arboricity of complete bipartite equipartition graphs was investigated in 2013. In this paper,
we study the exact value of the strong equitable vertex $3$-arboricity of
complete equipartition tripartite graphs.
\\[2mm]
\\[2mm]
\end{abstract}

\section{Introduction}

 In this paper, all graphs considered are finite and simple. For a
real number $x$, $\lceil x\rceil$ is the least integer not less than
$x$ and $\lfloor x \rfloor$ is the largest integer not larger than
$x$. We use $V(G)$, $E(G)$, $\delta(G)$ and $\Delta(G)$ to denote
the vertex set, edge set, minimum degree and maximum degree of $G$,
respectively. For a vertex $v \in V(G)$, let $ N_G(v)$ denote the
set of neighbors of $v$ in $G$ and $d_G(v)=|N_G(v)|$ denote the
degree of $v$. We use $d(v)$ instead of $d_G(v)$ for brevity. For
other undefined concepts, we refer the reader to \cite{Bollo}.

We associate positive integers $1,2,\ldots,t$ with colors and call
$f$ a \emph{$t$-coloring} of a graph $G$ if $f$ is a mapping from
$V(G)$ to $\{1,2,\ldots,t\}$. A t-coloring of $G$ is proper if any
two adjacent vertices have different colors. For $1 \leq i\leq t$,
let $V_i = \{v\,|\,f(v) = i\}$. A $t$-coloring of a graph $G$ is
said to be \emph{equitable} if $||V_i|-|V_j||\leq 1$ for all $i$ and
$j$, that is to say, every color class has size $\lfloor|V(G)|/t\rfloor$ or
$\lceil|V(G)|/t\rceil$. A graph $G$ is said to be \emph{properly
equitably $t$-colorable} if $G$ has a proper equitable $t$-coloring.
The smallest number $t$ for which $G$ is properly equitably
$t$-colorable is called \emph{the equitable chromatic number} of
$G$, denoted by $\chi^=(G)$.

The equitable coloring problem, introduced by Meyer
\cite{Meyer}, is motivated by a practical application to municipal
garbage collection \cite{Tucker}. In this context, the vertices of
the graph represent garbage collection routes. A pair of vertices
share an edge if the corresponding routes should not be run on the
same day. It is desirable that the number of routes ran on each day
be approximately the same. Therefore, the problem of assigning one
of the six weekly working days to each route reduces to finding a
proper equitable $6$-coloring. For more applications such as
scheduling, constructing timetables and load balance in parallel
memory systems, we refer to
\cite{Baker,Blazewicz,Das,Irani,Kitagawa,Smith}.

A properly equitably $t$-colorable graph may admit no proper equitable $t'$-colorings for some $t'> t$. For example, the complete bipartite graph $H := K_{2m+1,2m+1}$ has no proper equitable $(2m + 1)$-colorings, although it satisfies $\chi^=(H)=2$. This fact motivates another interesting parameter for proper equitable coloring. The equitable chromatic threshold of $G$, denoted by $\chi^กิ(G)$, is the smallest integer $t$ such that $G$ has a proper equitable $t'$-coloring for all $t'\geq t$. This notion was introduced by Fan et al. in \cite{Fan}.

In \cite{Fan}, Fan et al. considered relaxed equitable
coloring of graphs. They proved that every graph has a proper
equitable $\Delta(G)$-coloring such that each color class induces a
forest with maximum degree at most one. On the basis of this
research, Wu et al. \cite{Wu} introduced the notion of equitable
$(t, k)$-tree-coloring, which can be viewed as a generalization of
proper equitable $t$-coloring.

A \emph{$(t, k)$-tree-coloring} is a $t$-coloring $f$ of $G$
such that each component of $G[V_i]$ is a tree of maximum degree at
most $k$. An \emph{equitable $(t,
k)$-tree-coloring} is a $(t, k)$-tree-coloring that is equitable.
The \emph{equitable vertex $k$-arboricity} of a graph $G$, denoted
by ${va_k}^=(G)$, is the smallest integer $t$ such that $G$ has an
equitable $(t, k)$-tree-coloring. The \emph{strong equitable vertex
$k$-arboricity} of $G$, denoted by ${va_k}^\equiv(G)$, is the
smallest integer $t$ such that $G$ has an equitable $(t',
k)$-tree-coloring for every $t'\geq t$. It is clear that
${va_0}^=(G) = \chi^=(G)$ and ${va_0}^{\equiv}(G) =
{\chi}^{\equiv}(G) $ for every graph $G$, and ${va_k}^=(G)$ and
${va_k}^\equiv(G)$ may vary a lot.

In \cite{Wu}, Wu et al.
investigated the strong equitable vertex $k$-arboricity of complete
equipartition bipartite graphs and gave the bounds for ${va_1}^\equiv(K_{n,n})$ and ${va_\infty}^\equiv(K_{n,n})$. In this paper, we study the strong equitable vertex $3$-arboricity of complete equipartition tripartite graphs. In fact, we obtain the exact value of ${va_3}^\equiv(K_{n,n,n})$ for most cases.

\begin{thm}\label{th1}
If $n\equiv3 \ (mod~4)$, then ${va_3}^\equiv(K_{n,n,n})\leq3\lfloor
{\frac{n+1}{4}}\rfloor$.
\end{thm}

\begin{thm}\label{th2}

Let $k$ be a positive integer.

$(i)$ For $k\equiv1~(mod~5)$, ${va_3}^\equiv(K_{4k,4k,4k})=
\frac{12k+3}{5}$.

$(i i)$ For $k\equiv2~(mod~5)$, ${va_3}^\equiv(K_{4k,4k,4k})=
\frac{12k+6}{5}$.

$(iii)$ For $k\equiv3~(mod~5)$, ${va_3}^\equiv(K_{4k,4k,4k})=
\frac{12k+4}{5}$.

$(iv)$ For $k\equiv4~(mod~5)$, ${va_3}^\equiv(K_{4k,4k,4k})=
\frac{12k-3}{5}$.

$(v)$ For $k\equiv0~(mod~5)$, if $m\geq 2q+1$ where $q\geq 2$, then
${va_3}^\equiv(K_{4k,4k,4k})\leq 12m-3q-3$; \par if $m\geq 2q$ where
$q\geq 0$, then ${va_3}^\equiv(K_{4k,4k,4k})\leq 12m-3q$.
\end{thm}

\begin{thm}\label{th3}

Let $k$ be a positive integer.

$(i)$ For $k\equiv~2(mod~5)$, ${va_3}^\equiv(K_{4k+1,4k+1,4k+1})=
\frac{12k+6}{5}$.

$(i i)$ For $k\equiv3~(mod~5)$,
${va_3}^\equiv(K_{4k+1,4k+1,4k+1})= \frac{12k+9}{5}$.

$(iii)$ For $k\equiv4~(mod~5)$,
${va_3}^\equiv(K_{4k+1,4k+1,4k+1})= \frac{12k+7}{5}$.

$(iv)$ For $k\equiv0~(mod~5)$,
${va_3}^\equiv(K_{4k+1,4k+1,4k+1})= \frac{12k}{5}$.

$(v)$ For $k\equiv1~(mod~5)$, if $m\geq 2q+1$ where $q\geq 1$, then
${va_3}^\equiv(K_{4k+1,4k+1,4k+1})\leq$ \par $12m-3q$; if $m\geq 2q$
where $q\geq 0$, then ${va_3}^\equiv(K_{4k+1,4k+1,4k+1})\leq
12m-3q+3$.
\end{thm}

\begin{thm}\label{th4}

Let $k$ be a positive integer.

$(i)$ For $k\equiv3~(mod~5)$, ${va_3}^\equiv(K_{4k+2,4k+2,4k+2})=
\frac{12k+9}{5}$.

$(ii)$ For $k\equiv4~(mod~5)$,
${va_3}^\equiv(K_{4k+2,4k+2,4k+2})= \frac{12k+12}{5}$.

$(iii)$ For $k\equiv0~(mod~5)$,
${va_3}^\equiv(K_{4k+2,4k+2,4k+2})=\frac{12k+10}{5}$.

$(iv)$ For $k\equiv1~(mod~5)~(k\neq1)$,
${va_3}^\equiv(K_{4k+2,4k+2,4k+2})= \frac{12k+3}{5}$.

$(v)$ For $k\equiv2~(mod~5)$, if $m\geq 2q+1$ where $q\geq 1$, then
${va_3}^\equiv(K_{4k+2,4k+2,4k+2})\leq$ \par$12m-3q+3$; if $m\geq
2q$ where $q\geq 0$, then ${va_3}^\equiv(K_{4k+2,4k+2,4k+2})\leq
12m-3q+6$.
\end{thm}

\section{Preliminary}

\begin{pro}\label{pro2-1}
The complete $\ell$-partite graph $K_{n,\ldots,n}$ has an equitable
$(t, k)$-tree-coloring for every integer $t$ satisfying the
following condition:
$$
t={\ell}h \ (h\geq 1).
$$
\end{pro}

\begin{pf}
We can easily construct an equitable $(t, k)$-tree-coloring of $K_{n,\ldots,n}$ by dividing each partite set into $h$ classes equitably and coloring the vertices of each class with one color. \qed
\end{pf}

\begin{thm}\label{th2-2}
Let $K_{n,n,n}~(n\geq3)$ be a complete tripartite graph. Then
$$
{va_3}^\equiv(K_{n,n,n})\leq 3\left\lfloor
{\frac{n+1}{4}}\right\rfloor
$$
\end{thm}
\begin{pf}
By Proposition \ref{pro2-1}, in order to show
${va_3}^\equiv(K_{n,n,n})\leq 3\lfloor {\frac{n+1}{4}}\rfloor$, it
suffices to prove that $K_{n,n,n}$ has an equitable $(q,
3)$-tree-coloring for every $q$ satisfying the following condition:
$q\geq 3\lfloor {\frac{n+1}{4}}\rfloor +1$ and $3\nmid q$. Let
$X$,$Y$ and $Z$ be the partite sets of $K_{n,n,n}$.

Let $q=3a+1$ and $q\geq 3\lfloor {\frac{n+1}{4}}\rfloor +1$. Note
that $4q-3n\geq 12\lfloor {\frac{n+1}{4}}\rfloor +4-3n\geq 12
(\frac{n-2}{4})+4-3n=-2$.

If $q=\frac{3n-2}{4}$, then $4q-3n=-2$ and
$q=4+3t,~n=6+4t~(t=0,1,2,\cdots)$. Let $\{x_1, x_2, x_3, x_4, x_5,
x_6\} \subset X, \{y_1, y_2, y_3, y_4, y_5, y_6\} \subset Y$ and $
\{z_1, z_2, z_3, z_4, z_5, z_6\}\subset Z$. We color $x_1, x_2, x_3$
and $y_1$ with $1$ and color $y_2, y_3, y_4,y_5$ and $y_6$ with $2$
and color $x_4, x_5, x_6$ and $z_6$ with $3$ and color $z_1,
z_2,z_3,z_4 $ and $z_5$ with $4$. We divide each of
$X\backslash\{x_1, x_2, x_3, x_4, x_5, x_6\}$, $Y\backslash\{y_1,
y_2, y_3, y_4, y_5, y_6\}$ and $Z\backslash\{z_1, z_2, z_3, z_4,
z_5, z_6\}$ into $\frac{q-4}{3}$ classes equitably and color the
vertices of each class with a color in $\{5,\cdots,q\}$. Since
$\frac{n-6}{\frac{q-4}{3}}=4$, the resulting coloring is an
equitable $(q, 3)$-tree-coloring of $K_{n,n,n}$ with the size of
each color class being $4$ or $5$.

If $\frac{3n-2}{4} <q<\frac{3n+4}{4}$, then $\frac{3n-2}{4}
<3a+1<\frac{3n+4}{4}$, and hence $n-2<4a<n$. Then $4a=n-1$,
$a=\frac{n-1}{4}$ and $q=\frac{3n+1}{4}$. Therefore, $4q-3n=1$ and
$q=4+3t$, $n=5+4t~(t=0,1,2,\cdots)$. Let $\{x_1, x_2, x_3, x_4, x_5
\} \subset X$, $\{y_1, y_2, y_3, y_4, y_5\} \subset Y$ and $\{z_1,
z_2, z_3, z_4, z_5\}\subset Z$. We color $x_1, x_2, x_3$ and $y_1$
with $1$ and color $x_4, y_2, y_3$ and $y_4$ with $2$ and color
$x_5, z_1, z_2 $ and $z_3$ with $3$ and color $y_5, z_4 $ and $z_5$
with $4$. We divide each of $X\backslash\{x_1, x_2, x_3, x_4, x_5\},
Y\backslash\{y_1, y_2, y_3, y_4, y_5\}$ and $Z\backslash\{z_1, z_2,
z_3, z_4, z_5\}$ into $\frac{q-4}{3}$ classes equitably and color
the vertices of each class with a color in $\{5,\cdots,q\}$. Since
$\frac{n-5}{\frac{q-4}{3}}=4$, it follows that the resulting
coloring is an equitable $(q, 3)$-tree-coloring of $K_{n,n,n}$ with
the size of each color class being $3$ or $4$.

Suppose $\frac{3n+4}{4}\leq q\leq n$. Let $\{x_1, x_2, x_3, x_4 \}
\subset X, \{y_1, y_2, y_3, y_4\} \subset Y$ and $ \{z_1, z_2, z_3,
z_4\}\subset Z$. We color $x_1, x_2$ and $y_1$ with $1$ and color
$z_1, z_2$ and $y_2$ with $2$ and color $x_3, x_4 $ and $y_3$ with
$3$ and color $y_4, z_3 $ and $z_4$ with $4$. We divide each of
$X\backslash\{x_1, x_2, x_3, x_4\}, Y\backslash\{y_1, y_2, y_3,
y_4\}$ and $Z\backslash\{z_1, z_2, z_3, z_4\}$ into $\frac{q-4}{3}$
classes equitably and color the vertices of each class with a color
in $\{5,\cdots,q\}$. Since $\left\lfloor
{\frac{n-4}{\frac{q-4}{3}}}\right\rfloor=3$ and $\left\lceil
{\frac{n-4}{\frac{q-4}{3}}}\right\rceil=4$, it follows that the
resulting coloring is an equitable $(q, 3)$-tree-coloring of
$K_{n,n,n}$ with the size of each color class being $3$ or $4$.

Suppose $n+1\leq q\leq \frac{3n-1}{2}$. Let $e=xy$ be edge of
$K_{n,n,n}$ with $x\in X$ and $y\in Y$. We color $x, y$ with $1$ and
divide each of $X\backslash\{x\}$, $Y\backslash\{y\}$ and $Z$ into
$\frac{q-1}{3}$ classes equitably and color the vertices of each
class with a color in $\{2,\ldots,q\}$. One can easily check that
the resulting coloring is an equitable $(q,3)$-tree-coloring of
$K_{n,n,n}$ with the size of each color class being $2$ or $3$.

\noindent\textbf{Claim 1}: $q$ can not be in $(\frac{3n-1}{2},
\frac{3n+2}{2})$.

\noindent \textbf{Proof of Claim 1:} Assume, to the contrary, that
there exist a integer $q$ such that $q$ is in $(\frac{3n-1}{2},
\frac{3n+2}{2})$. Then
$$
\frac{3n-1}{2}< 3a+1 < \frac{3n+2}{2}.
$$
One can see that $n-1 < 2a < n$, a contradiction.\qed

Suppose $\frac{3n+2}{2}\leq q\leq 3(n-1)+1 $. Let $e=xy$ be an edge
of $K_{n,n,n}$ with $x\in X$ and $y\in Y$. We color $x$ and $y$ with
$1$ and divide each of $X\backslash\{x\}, Y\backslash\{y\}$ and $Z$
into $\frac{q-1}{3}$ classes equitably and color the vertices of
each class with a color in $\{2,\ldots,q\}$. One can easily check
that the resulting coloring is an equitable $(q, 3)$-tree-coloring
of $K_{n,n,n}$ with the size of each color class being $1$ or $2$.

Let $q=3a+2$ and $q\geq 3\lfloor {\frac{n+1}{4}}\rfloor+2$, we note
that $4q-3n\geq 12\lfloor {\frac{n+1}{4}}\rfloor +8-3n\geq 12
(\frac{n-2}{4})+8-3n=2$.

Suppose $\frac{3n+2}{4} \leq q\leq n$. Let $\{x_1,x_2\} \subset X,
\{y_1, y_2\} \subset Y$ and $ \{z_1, z_2\}\subset Z$. We color $x_1,
x_2$ and $y_1$ with $1$ and color $y_2, z_1$ and $z_2$ with $2$ and
divide each of $X\backslash\{x_1, x_2\}, Y\backslash\{y_1, y_2\}$
and $Z\backslash\{z_1, Z_2\}$ into $\frac{q-2}{3}$ classes equitably
and color the vertices of each class with a color in
$\{3,\ldots,q\}$. One can easily check that the resulting coloring
is an equitable $(q, 3)$-tree-coloring of $K_{n,n,n}$ with the size
of each color class being $3$ or $4$.

Suppose $n\leq q\leq \frac{3n-2}{2}$. Let $\{x_1, x_2\} \subset X,
\{y_1, y_2\} \subset Y$ and $ \{z_1, z_2\}\subset Z$. We color $x_1,
x_2$ and $y_1$ with $1$ and color $y_2, z_1$ and $z_2$ with $2$ and
divide each of $X\backslash\{x_1, x_2\}, Y\backslash\{y_1, y_2\}$
and $Z\backslash\{z_1, Z_2\}$ into $\frac{q-2}{3}$ classes equitably
and color the vertices of each class with a color in
$\{3,\ldots,q\}$. One can easily check that the resulting coloring
is an equitable $(q, 3)$-tree-coloring of $K_{n,n,n}$ with the size
of each color class being $2$ or $3$.

\noindent\textbf{Claim 2}: $q$ can no\emph{}t be in $(\frac{3n-2}{2},
\frac{3n+1}{2})$.

\noindent \textbf{Proof of Claim 2:} Assume, to the contrary, that
there exist a $q$ such that $q$ is in $(\frac{3n-2}{2},
\frac{3n+1}{2})$. Then,
$$
\frac{3n-2}{2}< 3a+2 < \frac{3n+1}{2}.
$$
We find that $n-2 < 2a < n-1$, a contradiction.\qed

Suppose $\frac{3n+1}{2}\leq q\leq 3(n-1)+2 $. Let $e=xy$ be an edge
of $K_{n,n,n}$ with $x\in X$ and $y\in Y$. Let $z$  be a vertex of
$Z$. We color $x$ and $y$ with $1$ and color $z$ with $2$   and
divide each of $X\backslash\{x\}, Y\backslash\{y\}$ and
$Z\backslash\{z\}$ into $\frac{q-2}{3}$ classes equitably and color
the vertices of each class with a color in $\{3,\ldots,q\}$. One can
easily check that the resulting coloring is an equitable $(q,
3)$-tree-coloring of $K_{n,n,n}$ with the size of each color class
being $1$ or $2$. \qed
\end{pf}

\section{Main results}

We now in a position to give our main results.

\subsection{The strong equitable vertex $3$-arboricity of $K_{4k+3,4k+3,4k+3}$}

From Theorem \ref{th2-2}, we can give a proof of Theorem \ref{th1}.

\noindent \textbf{Proof of Theorem $1$:} From Theorem \ref{th2-2},
we have ${va_3}^\equiv(K_{n,n,n})\leq 3\lfloor {\frac{n+1}{4}}\rfloor$.
\qed

Note that $K_{3,3,3}$ can attain the upper bound of Theorem $1$. One can check that ${va_3}^\equiv(K_{3,3,3})=3$.

\subsection{The strong equitable vertex $3$-arboricity of $K_{4k,4k,4k}$}

We investigate the strong equitable vertex $3$-arboricity of the
complete tripartite graph $K_{4k,4k,4k}$.

The following upper bound of $K_{4k,4k,4k}$ can be proved easily.

\begin{pro}\label{pro2}
If $k\geq 5m \ (m\geq 0)$ and $k\neq0$, then ${va_3}^\equiv(K_{4k,4k,4k})\leq
3k-3m$.
\end{pro}
\begin{pf}
We prove the proposition by induction on $m$. If $m=0$, then the result holds
by Theorem \ref{th2-2}. Assume $m\geq 1$. Since $k\geq 5m >5(m-1)$,
by the induction hypothesis and Proposition \ref{pro2-1}, we need to
prove that $K_{4k,4k,4k}$ has an equitable $(3k-3m+1,
3)$-tree-coloring and an equitable $(3k-3m+2, 3)$-tree-coloring.

Divide $X$ into $k-m+1$ classes equitably and color the vertices of
each class with a color in $\{1,2,\cdots,k-m+1\}$. Divide $Y$ into
$k-m$ classes equitably and color the vertices of each class with a
color in $\{k-m+2,\cdots,2k-2m+1\}$. Divide $Z$ into $k-m$ classes
equitably and color the vertices of each class with a color in
$\{2k-2m+2,\cdots,3k-3m+1\}$. It is easy to check that the resulting
coloring of $K_{4k,4k,4k}$ is an equitable $(3k-3m+1,
3)$-tree-coloring with the size of each color class being $4$ or
$5$.

Divide $X$ into $k-m+1$ classes equitably and color the vertices
of each class with a color in $\{1,2,\cdots,k-m+1\}$. Divide $Y$
into $k-m+1$ classes equitably and color the vertices of each class
with a color in $\{k-m+2,\cdots,2k-2m+2\}$. Divide $Z$ into $k-m$
classes equitably and color the vertices of each class with a color
in $\{2k-2m+3,\cdots,3k-3m+2\}$. It is easy to check that the
resulting coloring of $K_{4k,4k,4k}$ is an equitable $(3k-3m+2,
3)$-tree-coloring with the size of each color class being $4$ or
$5$. \qed
\end{pf}

The following corollaries are immediate.

\begin{cor}\label{cor1}
Let $k$ be a positive integer.

$(i )$ If $k\equiv1~(mod~5)$, then ${va_3}^\equiv(K_{4k,4k,4k})\leq
\frac{12k+3}{5}$.

$(ii)$ If $k\equiv2~(mod~5)$, then ${va_3}^\equiv(K_{4k,4k,4k})\leq
\frac{12k+6}{5}$.

$(iii)$ If $k\equiv3(mod~5)$, then ${va_3}^\equiv(K_{4k,4k,4k})\leq
\frac{12k+9}{5}$.

$(iv)$ If $k\equiv4~(mod~5)$, then ${va_3}^\equiv(K_{4k,4k,4k})\leq
\frac{12k+12}{5}$.
\end{cor}

\begin{lem}\label{lem1}
If $k\equiv 1(mod~5)$, then ${va_3}^\equiv(K_{4k,4k,4k})\geq
\frac{12k+3}{5}$.
\end{lem}

\begin{pf}
Let $k=5m+1$. We only need to show that $K_{4k,4k,4k}$ has no
equitable $(12m+2, 3)$-tree-coloring. Assume, to the contrary, that
$c$ is an equitable $(12m+2, 3)$-tree-coloring of $K_{4k,4k,4k}$.
Then the size of every color class in $c$ is at least $5$ because
$\lfloor\frac{12k}{12m+2}\rfloor=\lfloor\frac{60m+12}{12m+2}\rfloor=5$.

Let $c_i$ denote the number of those color
classes such that each color class contains exactly $i$ vertices, where $i=5,6$. Then we have the following two equations:
$$
5c_5+6c_6=60m+12
$$
$$
c_5+c_6=12m+2.
$$

We have the unique solution $c_5=12m$, $c_6=2$. Since each color class containing exactly $5$ or $6$ vertices must appear in some partite set of $K_{4k,4k,4k}$, it follows that $K_{4k,4k,4k}$ has no
equitable$(12m+2, 3)$-tree-coloring satisfying the above conditions. Then ${va_3}^\equiv(K_{4k,4k,4k})\geq\frac{12k+3}{5}$.\qed
\end{pf}

\begin{lem}\label{lem2}
If $k\equiv 2(mod~5)$, then ${va_3}^\equiv(K_{4k,4k,4k})\geq
\frac{12k+6}{5}$.
\end{lem}

\begin{pf}
Let $k=5m+2$. We only need to show that $K_{4k,4k,4k}$ has no
equitable $(12m+5, 3)$-tree-coloring. Assume, to the contrary, that
$c$ is an equitable $(12m+5, 3)$-tree-coloring of $K_{4k,4k,4k}$.
Then the size of every color class in $c$ is at least $4$ because
$\lfloor\frac{12k}{12m+5}\rfloor=\lfloor\frac{60m+24}{12m+5}\rfloor=4$.

Let $c_i$ denote the number of those color
classes such that each color class contains exactly $i$ vertices, where $i=4,5$. Then we have the following two equations:
$$
4c_4+5c_5=60m+24
$$
$$
c_4+c_5=12m+5.
$$

We have the unique solution $c_4=1$, $c_5=12m+4$. Since each color class containing exactly $5$ vertices must appear in some partite set of $K_{4k,4k,4k}$, it follows that $K_{4k,4k,4k}$ has no
equitable$(12m+5, 3)$-tree-coloring satisfying the above conditions. Then ${va_3}^\equiv(K_{4k,4k,4k})\geq\frac{12k+6}{5}$.\qed
\end{pf}

\begin{lem}\label{lem3}
If $k\equiv 3(mod~5)$, then ${va_3}^\equiv(K_{4k,4k,4k})\leq
\frac{12k+4}{5}$.
\end{lem}
\begin{pf}
Form $(iii)$ of Corollary\ref{cor1}, we have ${va_3}^\equiv(K_{4k,4k,4k})\leq
\frac{12k+9}{5}$. Let $k=5m+3$. We only need to show that $K_{4k,4k,4k}$ has an
equitable $(12m+8, 3)$-tree-coloring. Then the size of every color class is at least $4$ because $\lfloor\frac{12k}{12m+8}\rfloor=\lfloor\frac{60m+36}{12m+8}\rfloor=4$.

Let $c_i$ denote the number of those color
classes such that each color class contains exactly $i$ vertices, where $i=4,5$. Then we have the following two equations:
$$
4c_4+5c_5=60m+36
$$
$$
c_4+c_5=12m+8.
$$
We have the unique solution $c_4=4$, $c_5=12m+4$. Since each color class containing exactly $5$ vertices must appear in some partite set of $K_{4k,4k,4k}$, there are
 $4m+2$ color classes containing exactly $5$ vertices in some partite set of $K_{4k,4k,4k}$ and there are
 $4m+1$ color classes containing exactly $5$ vertices in other partite sets of $K_{4k,4k,4k}$. Since there are $20m+12$ vertices in every partite set of $K_{4k,4k,4k}$, there are $2$ vertices of color class containing exactly $4$ vertices in some partite set and there are $7$ vertices of color class containing exactly $4$ vertices in other partite sets. In this case, $K_{4k,4k,4k}$ has an equitable $(12m+8, 3)$-tree-coloring.\qed
\end{pf}

\begin{lem}\label{lem4}
If $k\equiv 3(mod~5)$, then ${va_3}^\equiv(K_{4k,4k,4k})\geq
\frac{12k+4}{5}$.
\end{lem}

\begin{pf}
Let $k=5m+3$. We only need to show that $K_{4k,4k,4k}$ has no
equitable $(12m+7, 3)$-tree-coloring. Assume, to the contrary, that
$c$ is an equitable $(12m+7, 3)$-tree-coloring of $K_{4k,4k,4k}$.
Then the size of every color class in $c$ is at least $5$ because
$\lfloor\frac{12k}{12m+7}\rfloor=\lfloor\frac{60m+36}{12m+7}\rfloor=5$.

Let $c_i$ denote the number of those color
classes such that each color class contains exactly $i$ vertices, where $i=5,6$. Then we have the following two equations:
$$
5c_5+6c_6=60m+36
$$
$$
c_5+c_6=12m+7.
$$

We have the unique solution $c_5=12m+6$, $c_6=1$. Since each color class containing exactly $5$ or $6$ vertices must appear in some partite set of $K_{4k,4k,4k}$, it follows that $K_{4k,4k,4k}$ has no
equitable$(12m+7, 3)$-tree-coloring satisfying the above conditions. Then ${va_3}^\equiv(K_{4k,4k,4k})\geq\frac{12k+4}{5}$.\qed
\end{pf}

\begin{lem}\label{lem5}
If $k\equiv 4(mod~5)$, then ${va_3}^\equiv(K_{4k,4k,4k})\leq
\frac{12k-3}{5}$.
\end{lem}

\begin{pf}
Form $(iv)$ of Corollary\ref{cor1}, we have ${va_3}^\equiv(K_{4k,4k,4k})\leq
\frac{12k+12}{5}$. Let $k=5m+4$. We need to show that $K_{4k,4k,4k}$ has an
equitable $(12m+11, 3)$-tree-coloring and an equitable $(12m+10, 3)$-tree-coloring.

If $K_{4k,4k,4k}$ has an equitable $(12m+11, 3)$-tree-coloring , then the size of every color class is at least $4$ because $\lfloor\frac{12k}{12m+11}\rfloor=\lfloor\frac{60m+48}{12m+11}\rfloor=4$.

Let $c_i$ denote the number of those color
classes such that each color class contains exactly $i$ vertices, where $i=4,5$. Then we have the following two equations:
$$
4c_4+5c_5=60m+48
$$
$$
c_4+c_5=12m+11.
$$
We have the unique solution $c_4=7$, $c_5=12m+4$. Since each color class containing exactly $5$ vertices must appear in some partite set of $K_{4k,4k,4k}$, there are
 $4m+2$ color classes containing exactly $5$ vertices in some partite set of $K_{4k,4k,4k}$ and there are
 $4m+1$ color classes containing exactly $5$ vertices in other partite sets of $K_{4k,4k,4k}$. Since there are $20m+16$ vertices in every partite set of $K_{4k,4k,4k}$, there are $6$ vertices of color class containing exactly $4$ vertices in some partite set and there are $11$ vertices of color class containing exactly $4$ vertices in other partite sets. In this case, $K_{4k,4k,4k}$ has an equitable $(12m+11, 3)$-tree-coloring.

 We can prove that  $K_{4k,4k,4k}$ has an equitable $(12m+10, 3)$-tree-coloring using an similar argument.

 Form the above argument and Proposition \ref{pro2-1}, we prove that ${va_3}^\equiv(K_{4k,4k,4k})\leq
\frac{12k-3}{5}$.\qed
\end{pf}

\begin{lem}\label{lem6}
If $k\equiv 4(mod~5)$, then ${va_3}^\equiv(K_{4k,4k,4k})\geq
\frac{12k-3}{5}$.
\end{lem}

\begin{pf}
Let $k=5m+4$. We only need to show that $K_{4k,4k,4k}$ has no
equitable $(12m+8, 3)$-tree-coloring. Assume, to the contrary, that
$c$ is an equitable $(12m+8, 3)$-tree-coloring of $K_{4k,4k,4k}$.
Then the size of every color class in $c$ is at least $5$ because
$\lfloor\frac{12k}{12m+8}\rfloor=\lfloor\frac{60m+48}{12m+8}\rfloor=5$.

Let $c_i$ denote the number of those color
classes such that each color class contains exactly $i$ vertices, where $i=5,6$. Then we have the following two equations:
$$
5c_5+6c_6=60m+48
$$
$$
c_5+c_6=12m+8.
$$

We have the unique solution $c_5=12m$, $c_6=8$. Since each color class containing exactly $5$ or $6$ vertices must appear in some partite set of $K_{4k,4k,4k}$, it follows that $K_{4k,4k,4k}$ has no
equitable$(12m+8, 3)$-tree-coloring satisfying the above conditions. Then ${va_3}^\equiv(K_{4k,4k,4k})\geq\frac{12k-3}{5}$.\qed
\end{pf}

In the following, we consider the remaining case $k\equiv~0(mod~5)$ and
give the proof of $(v)$ of Theorem \ref{th2}.

\begin{lem}\label{lem7}
If $k=5m$ and $m\geq 2$, then ${va_3}^\equiv(K_{4k,4k,4k})\leq
12m-3$.
\end{lem}
\begin{pf}
By Proposition \ref{pro2}, ${va_3}^\equiv(K_{4k,4k,4k})\leq 12m$. So
we need to prove that $K_{4k,4k,4k}$ has an equitable
$(12m-2,3)$-tree-coloring and an equitable $(12m-1,3)$-tree-coloring
by Proposition \ref{pro2-1}.

Divide $X$ into $4m-1$ classes equitably and color the vertices of
each class with a color in $\{1,2,\cdots,4m-1\}$. Divide $Y$ into
$4m-1$ classes equitably and color the vertices of each class with a
color in $\{4m,\cdots,8m-2\}$. Divide $Z$ into $4m$ classes
equitably and color the vertices of each class with a color in
$\{8m-1,\cdots,12m-2\}$. It is easy to check that the resulting
coloring of $K_{4k,4k,4k}$ is an equitable $(12m-2,
3)$-tree-coloring with the size of each color class being $5$ or
$6$.

Divide $X$ into $4m-1$ classes equitably and color the vertices of
each class with a color in $\{1,2,\cdots,4m-1\}$. Divide $Y$ into
$4m$ classes equitably and color the vertices of each class with a
color in $\{4m,\cdots,8m-1\}$. Divide $Z$ into $4m$ classes
equitably and color the vertices of each class with a color in
$\{8m,\cdots,12m-1\}$. It is easy to check that the resulting
coloring of $K_{4k,4k,4k}$ is an equitable $(12m-1,
3)$-tree-coloring with the size of each color class being $5$ or
$6$. \qed
\end{pf}

\begin{lem}\label{lem8}
If $k=5m$ and $m\geq 4$, then ${va_3}^\equiv(K_{4k,4k,4k})\leq
12m-6$.
\end{lem}
\begin{pf}
Similarly to the proof of Lemma
\ref{lem5}, we can prove it. \qed
\end{pf}

\begin{lem}\label{lem9}
If $k=5m$ and $m\geq 2q+1$ where $q\geq2$, then
${va_3}^\equiv(K_{4k,4k,4k})\leq 12m-3q-3$.
\end{lem}
\begin{pf}
We prove this theorem by induction on $q$. Suppose $q=2$. It is
equivalent to prove that if $m\geq 5$, then
${va_3}^\equiv(K_{4k,4k,4k})\leq 12m-9$. Since $m\geq 5>4$, it
follows Lemma \ref{lem6} from ${va_3}^\equiv(K_{4k,4k,4k})\leq
12m-6$. We need to prove that $K_{4k,4k,4k}$ has an equitable
$(12m-8,3)$-tree-coloring and an equitable $(12m-7,3)$-tree-coloring
by Proposition \ref{pro2-1}.

Divide $X$ into $4m-3$ classes equitably and color the vertices of
each class with a color in $\{1,2,\cdots,4m-3\}$. Divide $Y$ into
$4m-3$ classes equitably and color the vertices of each class with a
color in $\{4m-2,\cdots,8m-6\}$. Divide $Z$ into $4m-2$ classes
equitably and color the vertices of each class with a color in
$\{8m-5,\cdots,12m-8\}$. It is easy to check that the resulting
coloring of $K_{4k,4k,4k}$ is an equitable $(12m-8,
3)$-tree-coloring with the size of each color class being $5$ or
$6$.

Divide $X$ into $4m-3$ classes equitably and color the vertices of
each class with a color in $\{1,2,\cdots,4m-3\}$. Divide $Y$ into
$4m-2$ classes equitably and color the vertices of each class with a
color in $\{4m-2,\cdots,8m-5\}$. Divide $Z$ into $4m-2$ classes
equitably and color the vertices of each class with a color in
$\{8m-4,\cdots,12m-7\}$. It is easy to check that the resulting
coloring of $K_{4k,4k,4k}$ is an equitable $(12m-7,
3)$-tree-coloring with the size of each color class being $5$ or
$6$.

Suppose $q\geq 3$. Since $m\geq 2q+1>2(q-1)+1$, by the induction
hypothesis and Proposition \ref{pro2-1}, we need to prove that
$K_{4k,4k,4k}$ has an equitable $(12m-3q-1,3)$-tree-coloring and an
equitable $(12m-3q-2,3)$-tree-coloring.

Divide $X$ into $4m-q-1$ classes equitably and color the vertices of
each class with a color in $\{1,2,\cdots,4m-q-1\}$. Divide $Y$ into
$4m-q-1$ classes equitably and color the vertices of each class with
a color in $\{4m-q,\cdots,8m-2q-2\}$. Divide $Z$ into $4m-q$ classes
equitably and color the vertices of each class with a color in
$\{8m-2q-1,\cdots,12m-3q-2\}$. It is easy to check that the
resulting coloring of $K_{4k,4k,4k}$ is an equitable $(12m-3q-2,
3)$-tree-coloring with the size of each color class being $5$ or
$6$.

Divide $X$ into $4m-q-1$ classes equitably and color the vertices of
each class with a color in $\{1,2,\cdots,4m-q-1\}$. Divide $Y$ into
$4m-q$ classes equitably and color the vertices of each class with a
color in $\{4m-q,\cdots,8m-2q-1\}$. Divide $Z$ into $4m-q$ classes
equitably and color the vertices of each class with a color in
$\{8m-2q,\cdots,12m-3q-1\}$. It is easy to check that the resulting
coloring of $K_{4k,4k,4k}$ is an equitable $(12m-3q-1,
3)$-tree-coloring with the size of each color class being $5$ or
$6$. \qed
\end{pf}\\

\begin{lem}\label{lem10}
If $k=5m$ and $m\geq 2q$, where $q\geq0$, then
${va_3}^\equiv(K_{4k,4k,4k})\leq 12m-3q$.
\end{lem}
\begin{pf}
When $q=0,1,2$, the result holds by Proposition
\ref{pro2}, Lemma \ref{lem5} and Lemma \ref{lem6}. If $q\geq3$, then
$q-1\geq2$. Since $m\geq2q>2(q-1)+1$, it follows that
${va_3}^\equiv(K_{4k,4k,4k})\leq 12m-3(q-1)-3=12m-3q$. \qed
\end{pf}\\

\noindent \textbf{Proof of Theorem \ref{th2}:} Suppose $k\equiv1(mod~5)$.
From $(i)$ of Corollary \ref{cor1} and Lemma \ref{lem1},
${va_3}^\equiv(K_{4k,4k,4k})= \frac{12k+3}{5}$.
Suppose $k\equiv2(mod~5)$. From $(ii)$ of Corollary\ref{cor1} and Lemma \ref{lem2}, ${va_3}^\equiv(K_{4k,4k,4k})=\frac{12k+6}{5}$.
Suppose $k\equiv~3(mod~5)$. From Lemma \ref{lem3} and Lemma \ref{lem4}, ${va_3}^\equiv(K_{4k,4k,4k})=\frac{12k+4}{5}$.
Suppose $k\equiv~4(mod~5)$. From Lemma \ref{lem5} and Lemma \ref{lem6}, ${va_4}^\equiv(K_{4k,4k,4k})=\frac{12k-3}{5}$.
Suppose $k\equiv~0(mod~5)$. From Lemma \ref{lem9}, if $k=5m$ and $m\geq2q+1$ where $q\geq2$, then ${va_4}^\equiv(K_{4k,4k,4k})\leq 12m-3q-3$. From Lemma \ref{lem10}, if $k=5m$ and $m\geq 2q$ where $q\geq0$, then ${va_4}^\equiv(K_{4k,4k,4k})\leq 12m-3q$. \qed

\subsection{The strong equitable vertex $3$-arboricity of $K_{4k+1,4k+1,4k+1}$}

We investigate the strong equitable vertex $3$-arboricity of the
complete tripartite graph $K_{4k+1,4k+1,4k+1}$.

The following upper bound of $K_{4k+1,4k+1,4k+1}$ can be proved easily.

\begin{pro}\label{pro3}
If $k\geq 5m+1$ where $m\geq 0$, then ${va_3}^\equiv(K_{4k+1,4k+1,4k+1})\leq 3k-3m$.
\end{pro}
\begin{pf}
We prove the proposition by induction on $m$. If $m=0$, the result holds
by Theorem \ref{th2-2}.

Suppose $m\geq 1$. Since $k\geq 5m+1 >5(m-1)+1$, by the induction
hypothesis and Proposition \ref{pro2-1}, we need to prove that
$K_{4k+1,4k+1,4k+1}$ has an equitable $(3k-3m+1, 3)$-tree-coloring
and an equitable $(3k-3m+2, 3)$-tree-coloring.

Divide $X$ into $k-m+1$ classes equitably and color the vertices
of each class with a color in $\{1,2,\cdots,k-m+1\}$. Divide $Y$
into $k-m$ classes equitably and color the vertices of each class
with a color in $\{k-m+2,\cdots,2k-2m+1\}$. Divide $Z$ into $k-m$
classes equitably and color the vertices of each class with a color
in $\{2k-2m+2,\cdots,3k-3m+1\}$. It is easy to check that the
resulting coloring of $K_{4k+1,4k+1,4k+1}$ is an equitable
$(3k-3m+1, 3)$-tree-coloring with the size of each color class being
$4$ or $5$.

Divide $X$ into $k-m+1$ classes equitably and color the vertices
of each class with a color in $\{1,2,\cdots,k-m+1\}$. Divide $Y$
into $k-m+1$ classes equitably and color the vertices of each class
with a color in $\{k-m+2,\cdots,2k-2m+2\}$. Divide $Z$ into $k-m$
classes equitably and color the vertices of each class with a color
in $\{2k-2m+3,\cdots,3k-3m+2\}$. It is easy to check that the
resulting coloring of $K_{4k+1,4k+1,4k+1}$ is an equitable
$(3k-3m+2, 3)$-tree-coloring with the size of each color class being
$4$ or $5$. \qed
\end{pf}

The following corollaries are immediate.

\begin{cor}\label{cor2}

$(i)$ If $k\equiv~2(mod~5)$, then ${va_3}^\equiv(K_{4k+1,4k+1,4k+1})\leq  \frac{12k+6}{5}$.

$(ii )$ If $k\equiv~3(mod~5)$, then ${va_3}^\equiv(K_{4k+1,4k+1,4k+1})\leq \frac{12k+9}{5}$.

$(iii)$ If $k\equiv~4(mod~5)$, then ${va_3}^\equiv(K_{4k+1,4k+1,4k+1})\leq \frac{12k+12}{5}$.

$(iv)$ If $k\equiv~0(mod~5)$, then ${va_3}^\equiv(K_{4k+1,4k+1,4k+1})\leq \frac{12k+15}{5}$.
\end{cor}

\begin{lem}\label{lem11}
If $k\equiv 2(mod~5)$, then ${va_3}^\equiv(K_{4k+1,4k+1,4k+1})\geq
\frac{12k+6}{5}$.
\end{lem}

\begin{pf}
Let $k=5m+2$. We only need to show that $K_{4k+1,4k+1,4k+1}$ has no
equitable $(12m+5, 3)$-tree-coloring. Assume, to the contrary, that
$c$ is an equitable $(12m+5, 3)$-tree-coloring of $K_{4k+1,4k+1,4k+1}$.
Then the size of every color class in $c$ is exactly $5$ because
$\frac{12k+3}{12m+5}=\frac{60m+25}{12m+5}=5$. However, the size of every
partite sets is $20m+9$, which is not divisible by $5$, a contradiction. \qed
\end{pf}

\begin{lem}\label{lem12}
If $k\equiv 3(mod~5)$, then ${va_3}^\equiv(K_{4k+1,4k+1,4k+1})\geq
\frac{12k+9}{5}$.
\end{lem}

\begin{pf}
Let $k=5m+3$. We only need to show that $K_{4k+1,4k+1,4k+1}$ has no
equitable $(12m+8, 3)$-tree-coloring. Assume, to the contrary, that
$c$ is an equitable $(12m+8, 3)$-tree-coloring of $K_{4k+1,4k+1,4k+1}$.
Then the size of every color class in $c$ is at least $4$ because
$\lfloor\frac{12k+3}{12m+8}\rfloor=\lfloor\frac{60m+39}{12m+8}\rfloor=4$.

Let $c_i$ denote the number of those color
classes such that each color class contains exactly $i$ vertices, where $i=4,5$. Then we have the following two equations:
$$
4c_4+5c_5=60m+39
$$
$$
c_4+c_5=12m+8.
$$

We have the unique solution $c_4=1$, $c_5=12m+7$. Since each color class containing exactly $5$ vertices must appear in some partite set of $K_{4k+1,4k+1,4k+1}$, it follows that $K_{4k+1,4k+1,4k+1}$ has no
equitable$(12m+8, 3)$-tree-coloring satisfying the above conditions. Then ${va_3}^\equiv(K_{4k+1,4k+1,4k+1})\geq\frac{12k+9}{5}$.\qed
\end{pf}

\begin{lem}\label{lem13}
If $k\equiv 4(mod~5)$, then ${va_3}^\equiv(K_{4k+1,4k+1,4k+1})\leq
\frac{12k+7}{5}$.
\end{lem}
\begin{pf}
Form $(iii)$ of Corollary\ref{cor2}, we have ${va_3}^\equiv(K_{4k+1,4k+1,4k+1})\leq
\frac{12k+12}{5}$. Let $k=5m+4$. We only need to show that $K_{4k+1,4k+1,4k+1}$ has an equitable $(12m+11, 3)$-tree-coloring. Then the size of every color class is at least $4$ because $\lfloor\frac{12k+3}{12m+11}\rfloor=\lfloor\frac{60m+51}{12m+11}\rfloor=4$.

Let $c_i$ denote the number of those color
classes such that each color class contains exactly $i$ vertices, where $i=4,5$. Then we have the following two equations:
$$
4c_4+5c_5=60m+51
$$
$$
c_4+c_5=12m+11.
$$
We have the unique solution $c_4=4$, $c_5=12m+7$. Since each color class containing exactly $5$ vertices must appear in some partite set of $K_{4k+1,4k+1,4k+1}$, there are
 $4m+3$ color classes containing exactly $5$ vertices in some partite set of $K_{4k+1,4k+1,4k+1}$ and there are
 $4m+2$ color classes containing exactly $5$ vertices in other partite sets of $K_{4k+1,4k+1,4k+1}$. Since there are $20m+17$ vertices in every partite set of $K_{4k+1,4k+1,4k+1}$, there are $2$ vertices of color class containing exactly $4$ vertices in some partite set and there are $7$ vertices of color class containing exactly $4$ vertices in other partite sets. In this case, $K_{4k+1,4k+1,4k+1}$ has an equitable $(12m+11, 3)$-tree-coloring.\qed
\end{pf}

\begin{lem}\label{lem14}
If $k\equiv 4(mod~5)$, then ${va_3}^\equiv(K_{4k+1,4k+1,4k+1})\geq
\frac{12k+7}{5}$.
\end{lem}

\begin{pf}
Let $k=5m+4$. We only need to show that $K_{4k+1,4k+1,4k+1}$ has no
equitable $(12m+10, 3)$-tree-coloring. Assume, to the contrary, that
$c$ is an equitable $(12m+10, 3)$-tree-coloring of $K_{4k+1,4k+1,4k+1}$.
Then the size of every color class in $c$ is at least $5$ because
$\lfloor\frac{12k+3}{12m+10}\rfloor=\lfloor\frac{60m+51}{12m+10}\rfloor=5$.

Let $c_i$ denote the number of those color
classes such that each color class contains exactly $i$ vertices, where $i=5,6$. Then we have the following two equations:
$$
5c_5+6c_6=60m+51
$$
$$
c_5+c_6=12m+10.
$$

We have the unique solution $c_5=12m+9$, $c_6=1$. Since each color class containing exactly $5$ or $6$ vertices must appear in some partite set of $K_{4k+1,4k+1,4k+1}$, it follows that $K_{4k+1,4k+1,4k+1}$ has no
equitable$(12m+10, 3)$-tree-coloring satisfying the above conditions. Then ${va_3}^\equiv(K_{4k+1,4k+1,4k+1})\geq\frac{12k+7}{5}$.\qed
\end{pf}

\begin{lem}\label{lem15}
If $k\equiv 0(mod~5)$, then ${va_3}^\equiv(K_{4k+1,4k+1,4k+1})\leq
\frac{12k}{5}$.
\end{lem}

\begin{pf}
Form $(iv)$ of Corollary\ref{cor2}, we have ${va_3}^\equiv(K_{4k+1,4k+1,4k+1})\leq
\frac{12k+15}{5}$. Let $k=5m+5$. We need to show that $K_{4k+1,4k+1,4k+1}$ has an equitable $(12m+14, 3)$-tree-coloring and an equitable $(12m+13,3)$-tree-coloring.

If $K_{4k+1,4k+1,4k+1}$ has an equitable $(12m+14, 3)$-tree-coloring , then the size of every color class is at least $4$ because $\lfloor\frac{12k+3}{12m+14}\rfloor=\lfloor\frac{60m+63}{12m+14}\rfloor=4$.

Let $c_i$ denote the number of those color
classes such that each color class contains exactly $i$ vertices, where $i=4,5$. Then we have the following two equations:
$$
4c_4+5c_5=60m+63
$$
$$
c_4+c_5=12m+14.
$$
We have the unique solution $c_4=7$, $c_5=12m+7$. Since each color class containing exactly $5$ vertices must appear in some partite set of $K_{4k+1,4k+1,4k+1}$, there are
 $4m+3$ color classes containing exactly $5$ vertices in some partite set of $K_{4k+1,4k+1,4k+1}$ and there are
 $4m+2$ color classes containing exactly $5$ vertices in other partite sets of $K_{4k+1,4k+1,4k+1}$. Since there are $20m+21$ vertices in every partite set of $K_{4k+1,4k+1,4k+1}$, there are $6$ vertices of color class containing exactly $4$ vertices in some partite set and there are $11$ vertices of color class containing exactly $4$ vertices in other partite sets. In this case, $K_{4k+1,4k+1,4k+1}$ has an equitable $(12m+14, 3)$-tree-coloring.

 We can prove that  $K_{4k+1,4k+1,4k+1}$ has an equitable $(12m+13, 3)$-tree-coloring using an similar argument.

 Form the above argument and Proposition \ref{pro2-1}, we prove that ${va_3}^\equiv(K_{4k+1,4k+1,4k+1})\leq
\frac{12k}{5}$.\qed
\end{pf}

\begin{lem}\label{lem16}
If $k\equiv 0(mod~5)$, then ${va_3}^\equiv(K_{4k+1,4k+1,4k+1})\geq
\frac{12k}{5}$.
\end{lem}

\begin{pf}
Let $k=5m+5$. We only need to show that $K_{4k+1,4k+1,4k+1}$ has no
equitable $(12m+11, 3)$-tree-coloring. Assume, to the contrary, that
$c$ is an equitable $(12m+8, 3)$-tree-coloring of $K_{4k+1,4k+1,4k+1}$.
Then the size of every color class in $c$ is at least $5$ because
$\lfloor\frac{12k+3}{12m+11}\rfloor=\lfloor\frac{60m+63}{12m+11}\rfloor=5$.

Let $c_i$ denote the number of those color
classes such that each color class contains exactly $i$ vertices, where $i=5,6$. Then we have the following two equations:
$$
5c_5+6c_6=60m+63
$$
$$
c_5+c_6=12m+11.
$$

We have the unique solution $c_5=12m+3$, $c_6=8$. Since each color class containing exactly $5$ or $6$ vertices must appear in some partite set of $K_{4k+1,4k+1,4k+1}$, it follows that $K_{4k+1,4k+1,4k+1}$ has no
equitable$(12m+11, 3)$-tree-coloring satisfying the above conditions. Then ${va_3}^\equiv(K_{4k+1,4k+1,4k+1})\geq\frac{12k}{5}$.\qed
\end{pf}

In the following, we consider the remaining case $k\equiv~1(mod~5)$ and
give the proof of $(v)$ of Theorem \ref{th3}.

\begin{lem}\label{lem17}
If $k=5m+1$ and $m\geq2$, then ${va_3}^\equiv(K_{4k+1,4k+1,4k+1})\leq 12m$.
\end{lem}

\begin{pf}
By Proposition \ref{pro3} and Proposition \ref{pro2-1},
${va_3}^\equiv(K_{4k+1,4k+1,4k+1})\leq 12m+3$. It suffices to prove
that $K_{4k+1,4k+1,4k+1}$ has an equitable $(12m+2,3)$-tree-coloring
and an equitable $(12m+1,3)$-tree-coloring.

Divide $X$ into $4m+1$ classes equitably and color the vertices
of each class with a color in $\{1,2,\cdots,4m+1\}$. Divide $Y$
into $4m+1$ classes equitably and color the vertices of each class
with a color in $\{4m+2,\cdots,8m+2\}$. Divide $Z$ into $4m$
classes equitably and color the vertices of each class with a color
in $\{8m+3,\cdots,12m+2\}$. It is easy to check that the resulting
coloring of $K_{4k+1,4k+1,4k+1}$ is an equitable $(12m+2,
3)$-tree-coloring with the size of each color class being $5$ or
$6$.

Divide $X$ into $4m+1$ classes equitably and color the vertices
of each class with a color in $\{1,2,\cdots,4m+1\}$. Divide $Y$
into $4m$ classes equitably and color the vertices of each class
with a color in $\{4m+2,\cdots,8m+1\}$. Divide $Z$ into $4m$
classes equitably and color the vertices of each class with a color
in $\{8m+2,\cdots,12m+1\}$. It is easy to check that the resulting
coloring of $K_{4k+1,4k+1,4k+1}$ is an equitable $(12m+1,
3)$-tree-coloring with the size of each color class being $5$ or
$6$. \qed
\end{pf}

\begin{lem}\label{lem18}
If $k=5m+1$ and $m\geq 2q+1$ where $q\geq 1$, then
${va_3}^\equiv~(K_{4k+1,4k+1,4k+1})\leq 12m-3q$.
\end{lem}
\begin{pf}
We prove the theorem by induction on $q$. When $q=1$, it is
equivalent to prove that if $m\geq3$, then
${va_3}^\equiv(K_{4k,4k,4k})\leq 12m$. Since $m\geq3>2$, we have
${va_4}^\equiv(K_{4k+1,4k+1,4k+1})\leq 12m$ by Lemma \ref{lem13}. It
suffices to prove that $K_{4k+1,4k+1,4k+1}$ has an equitable
$(12m-1,3)$-tree-coloring and an equitable $(12m-2,3)$-tree-coloring
by Proposition \ref{pro2-1}.

Divide $X$ into $4m-1$ classes equitably and color the vertices
of each class with a color in $\{1,2,\cdots,4m-1\}$. Divide $Y$
into $4m$ classes equitably and color the vertices of each class
with a color in $\{4m,\cdots,8m-1\}$. Divide $Z$ into $4m$
classes equitably and color the vertices of each class with a color
in $\{8m,\cdots,12m-1\}$. It is easy to check that the resulting
coloring of $K_{4k+1,4k+1,4k+1}$ is an equitable $(12m-1,
3)$-tree-coloring with the size of each color class being $5$ or
$6$.

Divide $X$ into $4m-1$ classes equitably and color the vertices
of each class with a color in $\{1,2,\cdots,4m-1\}$. Divide $Y$
into $4m-1$ classes equitably and color the vertices of each class
with a color in $\{4m,\cdots,8m-2\}$. Divide $Z$ into $4m$
classes equitably and color the vertices of each class with a color
in $\{8m-1,\cdots,12m-2\}$. It is easy to check that the resulting
coloring of $K_{4k+1,4k+1,4k+1}$ is an equitable $(12m-2,
3)$-tree-coloring with the size of each color class being $5$ or
$6$. The result holds for $q=1$.

Suppose $q\geq 2$. Since $m\geq2q+1>2(q-1)+1$, by the induction
hypothesis and Proposition \ref{pro2-1}, we need to prove that
$K_{4k+1,4k+1,4k+1}$ has an equitable $(12m-3q+1,3)$-tree-coloring
and an equitable $(12m-3q+2,3)$-tree-coloring.

Divide $X$ into $4m-q+1$ classes equitably and color the vertices
of each class with a color in $\{1,2,\cdots,4m-q+1\}$. Divide $Y$
into $4m-q$ classes equitably and color the vertices of each class
with a color in $\{4m-q+2,\cdots,8m-2q+1\}$. Divide $Z$ into
$4m-q$ classes equitably and color the vertices of each class with a
color in $\{8m-2q+2,\cdots,12m-3q+1\}$. It is easy to check that the
resulting coloring of $K_{4k+1,4k+1,4k+1}$ is an equitable
$(12m-3q+1, 3)$-tree-coloring with the size of each color class
being $5$ or $6$.

Divide $X$ into $4m-q+1$ classes equitably and color the vertices
of each class with a color in $\{1,2,\cdots,4m-q+1\}$. Divide $Y$
into $4m-q+1$ classes equitably and color the vertices of each class
with a color in $\{4m-q+2,\cdots,8m-2q+2\}$. Divide $Z$ into
$4m-q$ classes equitably and color the vertices of each class with a
color in $\{8m-2q+3,\cdots,12m-3q+2\}$. It is easy to check that the
resulting coloring of $K_{4k+1,4k+1,4k=1}$ is an equitable
$(12m-3q+2, 3)$-tree-coloring with the size of each color class
being $5$ or $6$. \qed
\end{pf}

\begin{lem}\label{lem19}
If $k=5m+1$ and $m\geq2q$ where $q\geq0$, then
${va_3}^\equiv(K_{4k+1,4k+1,4k+1})\leq 12m-3q+3$.
\end{lem}
\begin{pf}
When $q=0,1$, the result holds by Proposition \ref{pro3} and Lemma
\ref{lem13}. If $q\geq 3$, then $q-1\geq 2$. Since
$m\geq2q>2(q-1)+1$, it follows that
${va_3}^\equiv(K_{4k+1,4k+1,4k+1})\leq 12m-3(q-1)=12m-3q+3$. \qed
\end{pf}\\

\noindent \textbf{Proof of Theorem \ref{th3}:} Suppose
$k\equiv2~(mod~5)$. From $(i)$ of Corollary \ref{cor2} and Lemma \ref{lem11},
${va_3}^\equiv(K_{4k+1,4k+1,4k+1})= \frac{12k+6}{5}$. Suppose
$k\equiv3~(mod~5)$. From $(ii)$ of Corollary \ref{cor2} and Lemma \ref{lem12},
${va_3}^\equiv(K_{4k+1,4k+1,4k+1})= \frac{12k+9}{5}$. Suppose
$k\equiv4~(mod~5)$. From Lemma \ref{lem13} and Lemma \ref{lem14},
${va_3}^\equiv(K_{4k+1,4k+1,4k+1})=\frac{12k+7}{5}$. Suppose
$k\equiv0~(mod~5)$. From Lemma \ref{lem15} and Lemma \ref{lem16},
${va_3}^\equiv(K_{4k+1,4k+1,4k+1})=\frac{12k}{5}$. Suppose
$k\equiv1~(mod~5)$. From Lemma \ref{lem18}, if $k=5m+1$ and $m\geq
2q+1$ where $q\geq1$, then ${va_3}^\equiv(K_{4k+1,4k+1,4k=1})\leq
12m-3q$. From Lemma \ref{lem19}, if $k=5m+1$ and $m\geq2q$ where
$q\geq0$, then ${va_3}^\equiv(K_{4k+1,4k+1,4k+1})\leq 12m-3q+3$.
\qed

\subsection{The strong equitable vertex $3$-arboricity of $K_{4k+2,4k+2,4k+2}$}

We investigate the strong equitable vertex $3$-arboricity of the
complete tripartite graph $K_{4k+2,4k+2,4k+2}$.

The following upper bound of $K_{4k+2,4k+2,4k+2}$ can be proved easily.

\begin{pro}\label{pro4}
If $k\geq 5m+2$ where $m\geq 0$, then ${va_3}^\equiv(K_{4k+2,4k+2,4k+2})\leq 3k-3m$.
\end{pro}

\begin{pf}
We prove the proposition  by induction on $m$. If $m=0$, the result holds
by Theorem \ref{th2-2}.

Suppose $m\geq 1$. Since $k\geq 5m+2 >5(m-1)+2$, by the induction
hypothesis and Proposition \ref{pro2-1}, we need to prove that
$K_{4k+2,4k+2,4k+2}$ has an equitable $(3k-3m+1,3)$-tree-coloring
and an equitable $(3k-3m+2,3)$-tree-coloring.

Divide $X$ into $k-m+1$ classes equitably and color the vertices of
each class with a color in $\{1,2,\cdots,k-m+1\}$. Divide $Y$ into
$k-m$ classes equitably and color the vertices of each class with a
color in $\{k-m+2,\cdots,2k-2m+1\}$. Divide $Z$ into $k-m$ classes
equitably and color the vertices of each class with a color in
$\{2k-2m+2,\cdots,3k-3m+1\}$. It is easy to check that the resulting
coloring of $K_{4k+2,4k+2,4k+2}$ is an equitable
$(3k-3m+1,3)$-tree-coloring with the size of each color class being
$4$ or $5$.

Divide $X$ into $k-m+1$ classes equitably and color the vertices of
each class with a color in $\{1,2,\cdots,k-m+1\}$. Divide $Y$ into
$k-m+1$ classes equitably and color the vertices of each class with
a color in $\{k-m+2,\cdots,2k-2m+2\}$. Divide $Z$ into $k-m$ classes
equitably and color the vertices of each class with a color in
$\{2k-2m+3,\cdots,3k-3m+2\}$. It is easy to check that the resulting
coloring of $K_{4k+2,4k+2,4k+2}$ is an equitable
$(3k-3m+2,3)$-tree-coloring with the size of each color class being
$4$ or $5$. \qed
\end{pf}

The following corollaries are immediate.

\begin{cor}\label{cor3}
$(i)$ If $k\equiv3~(mod~5)$, then
${va_3}^\equiv(K_{4k+2,4k+2,4k+2})\leq  \frac{12k+9}{5}$.

$(ii )$ If $k\equiv4~(mod~5)$, then
${va_3}^\equiv(K_{4k+2,4k+2,4k+2})\leq \frac{12k+12}{5}$.

$(iii)$ If $k\equiv0~(mod~5)$, then
${va_3}^\equiv(K_{4k+2,4k+2,4k+2})\leq \frac{12k+15}{5}$.

$(iv)$ If $k\equiv1~(mod~5)$, then
${va_3}^\equiv(K_{4k+2,4k+2,4k+2})\leq \frac{12k+18}{5}$.
\end{cor}

\begin{lem}\label{lem20}
If $k\equiv 3(mod~5)$, then ${va_3}^\equiv(K_{4k+2,4k+2,4k+2})\geq
\frac{12k+9}{5}$.
\end{lem}

\begin{pf}
Let $k=5m+3$. We only need to show that $K_{4k+2,4k+2,4k+2}$ has no
equitable $(12m+8, 3)$-tree-coloring. Assume, to the contrary, that
$c$ is an equitable $(12m+8, 3)$-tree-coloring of $K_{4k+2,4k+2,4k+2}$.
Then the size of every color class in $c$ is at least $5$ because
$\lfloor\frac{12k+6}{12m+8}\rfloor=\lfloor\frac{60m+42}{12m+8}\rfloor=5$.

Let $c_i$ denote the number of those color
classes such that each color class contains exactly $i$ vertices, where $i=5,6$. Then we have the following two equations:
$$
5c_5+6c_6=60m+42
$$
$$
c_5+c_6=12m+8.
$$

We have the unique solution $c_5=12m+6$, $c_6=2$. Since each color class containing exactly $5$ or $6$ vertices must appear in some partite set of $K_{4k+2,4k+2,4k+2}$, it follows that $K_{4k+2,4k+2,4k+2}$ has no
equitable$(12m+8, 3)$-tree-coloring satisfying the above conditions. Then ${va_3}^\equiv(K_{4k+2,4k+2,4k+2})\geq\frac{12k+9}{5}$.\qed
\end{pf}

\begin{lem}\label{lem21}
If $k\equiv 4(mod~5)$, then ${va_3}^\equiv(K_{4k+2,4k+2,4k+2})\geq
\frac{12k+12}{5}$.
\end{lem}

\begin{pf}
Let $k=5m+4$. We only need to show that $K_{4k+2,4k+2,4k+2}$ has no
equitable $(12m+11, 3)$-tree-coloring. Assume, to the contrary, that
$c$ is an equitable $(12m+11, 3)$-tree-coloring of $K_{4k+2,4k+2,4k+2}$.
Then the size of every color class in $c$ is at least $4$ because
$\lfloor\frac{12k+6}{12m+11}\rfloor=\lfloor\frac{60m+54}{12m+11}\rfloor=4$.

Let $c_i$ denote the number of those color
classes such that each color class contains exactly $i$ vertices, where $i=4,5$. Then we have the following two equations:
$$
4c_4+5c_5=60m+54
$$
$$
c_4+c_5=12m+11.
$$

We have the unique solution $c_4=1$, $c_5=12m+10$. Since each color class containing exactly $5$ vertices must appear in some partite set of $K_{4k+2,4k+2,4k+2}$, it follows that $K_{4k+2,4k+2,4k+2}$ has no
equitable$(12m+11, 3)$-tree-coloring satisfying the above conditions. Then ${va_3}^\equiv(K_{4k+2,4k+2,4k+2})\geq\frac{12k+12}{5}$.\qed
\end{pf}

\begin{lem}\label{lem22}
If $k\equiv 0(mod~5)$, then ${va_3}^\equiv(K_{4k+2,4k+2,4k+2})\leq
\frac{12k+10}{5}$.
\end{lem}
\begin{pf}
Form $(iii)$ of Corollary\ref{cor3}, we have ${va_3}^\equiv(K_{4k+2,4k+2,4k+2})\leq
\frac{12k+15}{5}$. Let $k=5m+5$. We only need to show that $K_{4k+2,4k+2,4k+2}$ has an
equitable $(12m+14, 3)$-tree-coloring. Then the size of every color class is at least $4$ because $\lfloor\frac{12k+6}{12m+14}\rfloor=\lfloor\frac{60m+66}{12m+14}\rfloor=4$.

Let $c_i$ denote the number of those color
classes such that each color class contains exactly $i$ vertices, where $i=4,5$. Then we have the following two equations:
$$
4c_4+5c_5=60m+66
$$
$$
c_4+c_5=12m+14.
$$
We have the unique solution $c_4=4$, $c_5=12m+10$. Since each color class containing exactly $5$ vertices must appear in some partite set of $K_{4k+2,4k+2,4k+2}$, there are
 $4m+4$ color classes containing exactly $5$ vertices in some partite set of $K_{4k+2,4k+2,4k+2}$ and there are
 $4m+3$ color classes containing exactly $5$ vertices in other partite sets of $K_{4k+2,4k+2,4k+2}$. Since there are $20m+22$ vertices in every partite set of $K_{4k+2,4k+2,4k+2}$, there are $2$ vertices of color class containing exactly $4$ vertices in some partite set and there are $7$ vertices of color class containing exactly $4$ vertices in other partite sets. In this case, $K_{4k+2,4k+2,4k+2}$ has an equitable $(12m+14, 3)$-tree-coloring.\qed
\end{pf}

\begin{lem}\label{lem23}
If $k\equiv 0(mod~5)$, then ${va_3}^\equiv(K_{4k+2,4k+2,4k+2})\geq
\frac{12k+10}{5}$.
\end{lem}

\begin{pf}
Let $k=5m+5$. We only need to show that $K_{4k+2,4k+2,4k+2}$ has no
equitable $(12m+13, 3)$-tree-coloring. Assume, to the contrary, that
$c$ is an equitable $(12m+13, 3)$-tree-coloring of $K_{4k+2,4k+2,4k+2}$.
Then the size of every color class in $c$ is at least $5$ because
$\lfloor\frac{12k+6}{12m+13}\rfloor=\lfloor\frac{60m+66}{12m+13}\rfloor=5$.

Let $c_i$ denote the number of those color
classes such that each color class contains exactly $i$ vertices, where $i=5,6$. Then we have the following two equations:
$$
5c_5+6c_6=60m+66
$$
$$
c_5+c_6=12m+13.
$$

We have the unique solution $c_5=12m+12$, $c_6=1$. Since each color class containing exactly $5$ or $6$ vertices must appear in some partite set of $K_{4k+2,4k+2,4k+2}$, it follows that $K_{4k+2,4k+2,4k+2}$ has no
equitable$(12m+13, 3)$-tree-coloring satisfying the above conditions. Then ${va_3}^\equiv(K_{4k+2,4k+2,4k+2})\geq\frac{12k+10}{5}$.\qed
\end{pf}

\begin{lem}\label{lem24}
If $k\equiv 1(mod~5)~(k\neq1)$, then ${va_3}^\equiv(K_{4k+2,4k+2,4k+2})\leq
\frac{12k+3}{5}$.
\end{lem}

\begin{pf}
Form $(iv)$ of Corollary\ref{cor3}, we have ${va_3}^\equiv(K_{4k+2,4k+2,4k+2})\leq
\frac{12k+18}{5}$. Let $k=5m+6$. We need to show that $K_{4k+2,4k+2,4k+2}$ has an
equitable $(12m+17, 3)$-tree-coloring and an equitable $(12m+16, 3)$-tree-coloring.

If $K_{4k+2,4k+2,4k+2}$ has an equitable $(12m+17, 3)$-tree-coloring , then the size of every color class is at least $4$ because $\lfloor\frac{12k+6}{12m+17}\rfloor=\lfloor\frac{60m+78}{12m+17}\rfloor=4$.

Let $c_i$ denote the number of those color
classes such that each color class contains exactly $i$ vertices, where $i=4,5$. Then we have the following two equations:
$$
4c_4+5c_5=60m+78
$$
$$
c_4+c_5=12m+17.
$$
We have the unique solution $c_4=7$, $c_5=12m+10$. Since each color class containing exactly $5$ vertices must appear in some partite set of $K_{4k+2,4k+2,4k+2}$, there are
 $4m+4$ color classes containing exactly $5$ vertices in some partite set of $K_{4k+2,4k+2,4k+2}$ and there are
 $4m+3$ color classes containing exactly $5$ vertices in other partite sets of $K_{4k+2,4k+2,4k+2}$. Since there are $20m+26$ vertices in every partite set of $K_{4k+2,4k+2,4k+2}$, there are $6$ vertices of color class containing exactly $4$ vertices in some partite set and there are $11$ vertices of color class containing exactly $4$ vertices in other partite sets. In this case, $K_{4k+2,4k+2,4k+2}$ has an equitable $(12m+17, 3)$-tree-coloring.

 We can prove that  $K_{4k+2,4k+2,4k+2}$ has an equitable $(12m+16, 3)$-tree-coloring using an similar argument.

 Form the above argument and Proposition \ref{pro2-1}, we prove that ${va_3}^\equiv(K_{4k+2,4k+2,4k+2})\leq
\frac{12k+3}{5}$.\qed
\end{pf}

\begin{lem}\label{lem25}
If $k\equiv 1(mod~5~(k\neq1)$, then ${va_3}^\equiv(K_{4k+2,4k+2,4k+2})\geq
\frac{12k+3}{5}$.
\end{lem}

\begin{pf}
Let $k=5m+6$. We only need to show that $K_{4k+2,4k+2,4k+2}$ has no
equitable $(12m+14, 3)$-tree-coloring. Assume, to the contrary, that
$c$ is an equitable $(12m+14, 3)$-tree-coloring of $K_{4k+2,4k+2,4k+2}$.
Then the size of every color class in $c$ is at least $5$ because
$\lfloor\frac{12k+6}{12m+14}\rfloor=\lfloor\frac{60m+78}{12m+14}\rfloor=5$.

Let $c_i$ denote the number of those color
classes such that each color class contains exactly $i$ vertices, where $i=5,6$. Then we have the following two equations:
$$
5c_5+6c_6=60m+78
$$
$$
c_5+c_6=12m+14.
$$

We have the unique solution $c_5=12m+6$, $c_6=8$. Since each color class containing exactly $5$ or $6$ vertices must appear in some partite set of $K_{4k+2,4k+2,4k+2}$, it follows that $K_{4k+2,4k+2,4k+2}$ has no
equitable$(12m+14, 3)$-tree-coloring satisfying the above conditions. Then ${va_3}^\equiv(K_{4k+2,4k+2,4k+2})\geq\frac{12k+3}{5}$.\qed
\end{pf}

In the following, we consider the remaining case $k\equiv~2(mod~5)$ and
give the proof of $(v)$ of Theorem \ref{th4}.

\begin{lem}\label{lem26}
If $k=5m+2$ and $m\geq2$, then ${va_3}^\equiv(K_{4k+2,4k+2,4k+2})\leq 12m+3$.
\end{lem}

\begin{pf}
By Proposition \ref{pro4} and Proposition \ref{pro2-1},
${va_3}^\equiv(K_{4k+2,4k+2,4k+2})\leq 12m+6$. So we need to prove that
$K_{4k+2,4k+2,4k+2}$ has an equitable $(12m+4,3)$-tree-coloring and
an equitable $(12m+5,3)$-tree-coloring.

Divide $X$ into $4m+1$ classes equitably and color the vertices
of each class with a color in $\{1,2,\cdots,4m+1\}$. Divide $Y$
into $4m+1$ classes equitably and color the vertices of each class
with a color in $\{4m+2,\cdots,8m+2\}$. Divide $Z$ into $4m+2$
classes equitably and color the vertices of each class with a color
in $\{8m+3,\cdots,12m+4\}$. It is easy to check that the resulting
coloring of $K_{4k+2,4k+2,4k+2}$ is an equitable $(12m+4,
3)$-tree-coloring with the size of each color class being $5$ or
$6$.

Divide $X$ into $4m+1$ classes equitably and color the vertices
of each class with a color in $\{1,2,\cdots,4m+1\}$. Divide $Y$
into $4m+2$ classes equitably and color the vertices of each class
with a color in $\{4m+2,\cdots,8m+3\}$. Divide $Z$ into $4m+2$
classes equitably and color the vertices of each class with a color
in $\{8m+4,\cdots,12m+5\}$. It is easy to check that the resulting
coloring of $K_{4k+2,4k+2,4k+2}$ is an equitable $(12m+5,
3)$-tree-coloring with the size of each color class being $5$ or
$6$. \qed
\end{pf}

\begin{lem}\label{lem27}
If $k=5m+2$ and $m\geq2q+1$ where $q\geq 1$, then
${va_3}^\equiv(K_{4k+2,4k+2,4k+2})\leq 12m-3q+3$.
\end{lem}

\begin{pf}
We prove this theorem by induction on $q$. When $q=1$, it is equivalent
to prove that if $m\geq3$, then
${va_3}^\equiv(K_{4k+2,4k+2,4k+2})\leq 12m$. Since
$m\geq3>2$, it follows from Lemma
\ref{lem20} that ${va_3}^\equiv(K_{4k+2,4k+2,4k+2})\leq 12m+3$. We need to prove that $K_{4k+2,4k+2,4k+2}$ has an
equitable $(12m+1,3)$-tree-coloring and an equitable
$(12m+2,3)$-tree-coloring by Proposition \ref{pro2-1}.

Divide $X$ into $4m+1$ classes equitably and color the vertices
of each class with a color in $\{1,2,\cdots,4m+1\}$. Divide $Y$
into $4m$ classes equitably and color the vertices of each class
with a color in $\{4m+2,\cdots,8m+1\}$. Divide $Z$ into $4m$
classes equitably and color the vertices of each class with a color
in $\{8m+2,\cdots,12m+1\}$. It is easy to check that the resulting
coloring of $K_{4k+1,4k+1,4k+1}$ is an equitable $(12m+1,
3)$-tree-coloring with the size of each color class being $5$ or
$6$.

Divide $X$ into $4m+1$ classes equitably and color the vertices
of each class with a color in $\{1,2,\cdots,4m+1\}$. Divide $Y$
into $4m+1$ classes equitably and color the vertices of each class
with a color in $\{4m+2,\cdots,8m+2\}$. Divide $Z$ into $4m$
classes equitably and color the vertices of each class with a color
in $\{8m+3,\cdots,12m+2\}$. It is easy to check that the resulting
coloring of $K_{4k+2,4k+2,4k+2}$ is an equitable $(12m+2,
3)$-tree-coloring with the size of each color class being $5$ or
$6$.

Suppose $q\geq2$. Since $m\geq2q+1>2(q-1)+1$, by the induction
hypothesis and Proposition \ref{pro2-1}, it suffices to prove that
$K_{4k+2,4k+2,4k+2}$ has an equitable $(12m-3q+5,3)$-tree-coloring
and an equitable $(12m-3q+4,3)$-tree-coloring.

Divide $X$ into $4m-q+1$ classes equitably and color the vertices
of each class with a color in $\{1,2,\cdots,4m-q+1\}$. Divide $Y$
into $4m-q+1$ classes equitably and color the vertices of each class
with a color in $\{4m-q+2,\cdots,8m-2q+2\}$. Divide $Z$ into
$4m-q+2$ classes equitably and color the vertices of each class with
a color in $\{8m-2q+3,\cdots,12m-3q+4\}$. It is easy to check that
the resulting coloring of $K_{4k+2,4k+2,4k+2}$ is an equitable
$(12m-3q+4,3)$-tree-coloring with the size of each color class
being $5$ or $6$.

Divide $X$ into $4m-q+1$ classes equitably and color the vertices
of each class with a color in $\{1,2,\cdots,4m-q+1\}$. Divide $Y$
into $4m-q+2$ classes equitably and color the vertices of each class
with a color in $\{4m-q+2,\cdots,8m-2q+3\}$. Divide $Z$ into
$4m-q+2$ classes equitably and color the vertices of each class with
a color in $\{8m-2q+4,\cdots,12m-3q+5\}$. It is easy to check that
the resulting coloring of $K_{4k+2,4k+2,4k+2}$ is an equitable
$(12m-3q+5,3)$-tree-coloring with the size of each color class
being $5$ or $6$. \qed
\end{pf}
\begin{lem}\label{lem28}
If $k=5m+2$ and $m\geq2q$ where $q\geq0$, then ${va_3}^\equiv(K_{4k+2,4k+2,4k+2})\leq 12m-3q+6$.
\end{lem}
\begin{pf}
When $q=0,1$, the result holds by Proposition \ref{pro4} and
Lemma \ref{lem20}. If $q\geq2$, then $q-1\geq1$. Since
$m\geq2q>2(q-1)+1$, it follows that ${va_3}^\equiv(K_{4k+2,4k+2,4k+2})\leq
12m-3(q-1)+3=12m-3q+6$. \qed
\end{pf}\\

\noindent \textbf{Proof of Theorem \ref{th4}:} Suppose $k\equiv3(mod 5)$.
From $(i)$ of Corollary \ref{cor3} and Lemma \ref{lem20},
${va_3}^\equiv(K_{4k+2,4k+2,4k+2})=\frac{12k+9}{5}$.
Suppose $k\equiv4(mod 5)$. From $(ii)$ of
Corollary\ref{cor3} and Lemma \ref{lem21}, ${va_3}^\equiv(K_{4k+2,4k+2,4k+2})= \frac{12k+12}{5}$.
Suppose $k\equiv0(mod 5)$. From Lemma \ref{lem22} and Lemma \ref{lem23}, ${va_3}^\equiv(K_{4k+2,4k+2,4k+2})= \frac{12k+10}{5}$.
Suppose $k\equiv1(mod 5)$. From Lemma \ref{lem24} and Lemma \ref{lem25}, ${va_3}^\equiv(K_{4k+2,4k+2,4k+2})= \frac{12k+5}{5}$.
Suppose $k\equiv2(mod~5)$. From Lemma \ref{lem27}, if $k=5m+2$ and $m\geq2q+1$, where $q\geq 1$, then ${va_3}^\equiv(K_{4k+2,4k+2,4k+2})\leq 12m-3q+3$.
From Lemma \ref{lem28}, if $k=5m+2$ and $m\geq2q$, where $q\geq0$, then ${va_3}^\equiv(K_{4k+2,4k+2,4k+2})\leq 12m-3q+6$. \qed

\end{document}